\input amssym.def
\input amssym

\input prepictex
\input pictex
\input postpictex

\def\sq{\hbox{\rlap{$\sqcap$}$\sqcup$}}
\def\qed{\ifmmode\sq\else{\unskip\nobreak\hfil
\penalty50\hskip1em\null\nobreak\hfil\sq
\parfillskip=0pt\finalhyphendemerits=0\endgraf}\fi}

\def\bbbc{ I \!\!\!\! C}
\def\nz{ I \! N}

\def\pr{ I \!\! P}
\def\r{\rightarrow}

\def\p{\partial}
\def\s{\subset}

\def\o{\overline}

\def\D{\Delta}
\def\al{\alpha}
\def\be{\beta}

\def\se{\setminus}
\def\de{\delta}

\def\la{\lambda}
\def\ep{\epsilon}
\def\ga{\gamma}
\def\la{\lambda}

\def\Ga{\Gamma}
\def\v{\vartheta}

\catcode`\"=\active
\def"{\"}
\def\3{\char\ss}

\def\bbbc{\Bbb C}
\def\nz{\Bbb N}

\def\pr{\Bbb P}

\documentstyle[12pt]{article}

\newtheorem{defi}{Definition}[section]
\newtheorem{prop}[defi]{Proposition}
\newtheorem{theo}[defi]{Theorem}
\newtheorem{lem}[defi]{Lemma}

\newtheorem{cor}[defi]{Corollary}
\newtheorem{rem}[defi]{Remark}
\newtheorem{conj}[defi]{Conjecture}
\begin{document}

\begin{center}
{\LARGE Hyperbolicity of the Complements}\\
{\LARGE of Plane \vspace{1cm} Algebraic Curves}

{\large G.Dethloff, G.Schumacher, P.M.Wong}
\end{center}
\tableofcontents

\section{Introduction}

Hyperbolic manifolds have been studied in complex analysis as the
generalizations of hyperbolic Riemann surfaces to higher dimensions.
Moreover, the theory of hyperbolic manifolds is closely related to other
areas (cf.\ eg. \cite{LA1}).

However, only very few quasi-projective (non closed) hyperbolic manifolds
are known.  But one still believes that e.g.\ the complements of `most'
hypersurfaces in $\pr_n$ are hyperbolic, if only their degree is at
least 2n+1, more precisely:

\begin{conj}
Let ${\cal C}(d_1,\ldots ,d_k)$ be the space of $k$ tupels of
hypersurfaces $\,\Ga = (\Ga_1 , \ldots , \Ga_k )\,$ in $\pr_n$, where
${\rm deg}(\Ga_i)=d_i$.  Then for all $(d_1,\ldots ,d_k)$ with $\,
\sum_{i=1}^k d_i =:d \geq 2n+1\,$ the set $\,{\cal H}(d_1,\ldots ,d_k)=
\{ \Ga \in {\cal C}(d_1,\ldots ,d_k) :  \pr_n \se \bigcup_{i=1}^k
\Ga_i\, $ {\rm is complete hyperbolic and hyperbolically embedded}$\}\,$
contains the complement of a proper algebraic subset of ${\cal
C}(d_1,\ldots ,d_k)$.
\end{conj}

For complements of hypersurfaces in $\pr_n$ this was posed by Kobayashi
as `Problem 3' in his book \cite{KO}, and later by Zaidenberg in his
paper \cite{ZA}.

In this paper, we shall deal with the complements of plane curves i.e.
the case n=2.

Other than in the case of 5 lines $({\cal C}(1,1,1,1,1))$, the
conjecture was previously proved by M.~Green in \cite{GRE2} in the case
of a curve $\Ga$ consisting of one quadric and three lines (${\cal
C}(2,1,1,1)$).  Furthermore, it was shown for ${\cal
C}(d_1,\ldots,d_k)$, whenever $k\geq5$, by Babets in \cite{BA}.  A
closely related result by Green in \cite{GRE1} is that for any four
non-redundant hypersurfaces $\Gamma_j$, $j=1,\ldots 4$ in $\pr_2$ any
entire curve $f:\bbbc \to \pr_2 \setminus \bigcup_{j=1}^4 \Gamma_j$ is
algebraically degenerate.  (The degeneracy locus of the Kobayashi
pseudometric was studied by Adachi and Suzuki in \cite{A--S1},
\cite{A--S2}).

In fact, for generic configurations, any such algebraically degenerate
map is constant, hence the conjecture is true for any family ${\cal
C}(d_1,\ldots,\linebreak d_k)$ with $k \geq 4$ (cf.\ cf.\  Theorem \ref{4c}).
This includes the case of a curve $\Ga$ consisting of 2 quadrics and 2
lines.  We also give another proof of Green's result, which yields a
slightly stronger result related to the statement of a second main
theorem of value distribution theory in this situation.

It seems that the conjecture is the more difficult the smaller k is.
Already the case k=3 seems to be very hard:  In 1989 H.~Grauert worked
on the case of a curve $\Ga$ consisting of 3 quadrics, i.e.  ${\cal
C}(2,2,2)$, in \cite{GR}, using sophisticated differential geometric
methods including Jet-metrics.  We believe that the methods developed
there might be suited for proving major parts of the conjecture.  For the
time being, however, certain technical problems still exist with these
methods including the case ${\cal C}(2,2,2)$.

The main result of this paper (Theorem \ref{mt}) is a proof of the
conjecture for 3 quadrics.  Our
methods are completely different from those used in \cite{GR} --- instead of
differential geometry we use value distribution theory:
For any pair of quadrics which intersect transversally,
there are 6 lines through the intersection points, out of which 4
are in general position. We first show that we can assign a set of 12
lines in general position to any generic system $C$ of 3 quadrics.
Let $f: \bbbc \r \pr_2  \se C$ be an entire holomorphic curve.
Our method now essentially consists of showing that the defect of $f$
with respect to the above 12 lines had to be at least equal to 4
unless $f$ is algebraically degenerate. (For technical reasons our
exposition is based on the Second Main Theorem rather than the
defect relation). The last step is to show that this fact is
actually sufficient for generic complements of 3 quadrics to be
complete hyperbolic and hyperbolically embedded.

For ${\cal C}(2,2,1)$, i.e. two quadrics and a line, our result
states the existence of an open set, which contains a quasi-projective
set of codimension one, of configurations, where the conjecture is true
(Theorem $\ref{thm221}$).
The somewhat lengthy proof is based on a generalized Borel lemma.
With the same methods we prove that also the
complement of three generic Fermat quadrics is hyperbolic.

The paper is organized as follows:  In section~2 we collect, for the
convenience of the reader, some basics from value distribution theory,
and, in section~3, some consequences from Brody's techniques for later
reference.  In section~4 we prove some `algebraic' hyperbolicity of
generic complements of certain curves.  Next, in section~5 we prove Theorem
\ref{4c}.  In section~6
we study linear systems of lines associated to systems of 3 quadrics.
Section~7 contains the proof of Theorem~\ref{mt}.  In
section~8 we treat complements of two quadrics and a line and
complements of three Fermat quadrics.

The first named author would like to thank S.~Frankel (Nantes),
H.~Grauert (G"ottingen), S.~Kosarew (Grenoble) and M.~Zaidenberg
(Grenoble) for valuable discussions, the Department of Mathematics at
Notre Dame for its hospitality, and the DFG, especially the `Schwer\-punkt
Kom\-ple\-xe Man\-nig\-fal\-tig\-kei\-ten' in Bochum for support.  The second
named author would like to thank H.~Grauert, W.Stoll (Notre Dame) and
M.Zaidenberg for valuable discussions, and the Department of Mathematics
at Notre Dame and the SFB~170 in G\"{o}ttingen for its hospitality and
the Schwerpunkt `Komplexe Mannigfaltig\-kei\-ten' for support.  The third
named author would like to thank the SFB 170 and the NSF for partial
support.


\section{Some tools from Value Distribution Theory}

In this section we fix some notations and quote some facts from Value
Distribution Theory.  We give references but do not trace these facts
back to the original papers.

We define the characteristic function and the counting function, and
give some formulas for these.

Let $\,||z||^2= \sum_{j=0}^n |z_j|^2$, where $(z_0,\ldots ,z_n) \in
\bbbc^{n+1}$, let $\D_t = \{\xi \in \bbbc :  |\xi| < t \}$, and let $d^c =
(i/4 \pi) (\overline{\p} - \p)$.  Let $r_0$ be a fixed positive
number and let $\,r \geq r_0$.  Let $\,f:\bbbc \r \pr_n\,$ be entire, i.e.
$f$ can be written as $\, f=[f_0:\ldots :f_n]\,$ with holomorphic
functions $\, f_j :  \bbbc \r \bbbc\, , j=0,\ldots ,n\,$ without common
zeroes.  Then the {\it characteristic function} $T(f,r)$ is defined as
$$ T(f,r) = \int_{r_0}^r \frac{dt}{t} \int_{\D_t} dd^c \log ||f||^2$$
Let furthermore $\, D=\{ P=0\}\,$ be a divisor in $\pr_n$, given by a
homogeneous polynomial $P$.  Assume $\, f(\bbbc) \not\s \hbox{ {\rm
support}}(D)$.  Let $\,n_f(D,t)\,$ denote the number of zeroes of $\, P
\circ f\,$ inside $\, \D_t\,$ (counted with multiplicities).  Then we
define the {\it counting function} as
$$
N_f(D,r) = \int_{r_0}^r n_f(D,t) \frac{dt}{t}
$$

Stokes Theorem and transformation to polar coordinates imply (cf.\ \cite{WO}):
\begin{equation} \label{1}
T(f,r) =
\frac{1}{4 \pi} \int_0^{2 \pi} \log  ||f||^2 (re^{i \v})d \v + O(1).
\end{equation}

The characteristic function as defined by Nevanlinna for a holomorphic
function $\,f:  \bbbc \r \bbbc$ is
$$
T_0(f,r) = \frac{1}{2 \pi} \int_0^{2 \pi}
\log ^+ |f(re^{i \v})| d \v .
$$
For the associated map  $\, [f:1]:  \bbbc \r \pr_1$ one has
\begin{equation} \label{2}
T_0(f,r) = T([1:f],r) + O(1)
\end{equation}
(cf.\ \cite{HA}).

By abuse of notation we will, from now on, for a function $\, f:  \bbbc \r
\bbbc$, write $T(f,r)$ instead of $T_0(f,r)$.  Furthermore we
sometimes use $N(f,r)$ instead of $N_f([z_0=0],r)$.

The concept of finite order is essential for later applications.
\begin{defi}
Let $s(r)$ be a positive, monotonically increasing function
defined for $\,r \geq r_0$. If
$$ \overline{\lim_{r \r \infty} } \frac{\log  s(r)}{\log  r} = \la$$
then $s(r)$ is said to be of order $\la$. For entire $\,f:\bbbc \r \pr_n\,$
or $\, f: \bbbc \r \bbbc\,$ we say that $f$ is of order $\la$, if
$T(f,r)$ is.
\end{defi}

\begin{rem}\label{remfo}
Let $f=[f_0:\ldots:f_n]:\bbbc \to \pr_n$ be a
holomorphic map of finite order $\la$.  Then $\log T(f,r)= O(\log r)$.
\end{rem}

We need the following:
\begin{lem} \label{e}
Assume that $\,f: \bbbc \r \pr_n\,$ is an entire map and misses the
divisors
$\,\{ z_j = 0\}\,$ for $j=0,\ldots,n$ (i.e. the coordinate hyperplanes
of $\pr_n$).
Assume that $f$ has order at most $\la$. Then $f$ can be written as
$\,f = [1:f_1:\ldots :f_n]\,$ with $\, f_j(\xi) = e^{P_j(\xi)}$, where
the $P_j(\xi)$
are polynomials in $\xi$ of degree $d_j\leq \lambda$.
\end{lem}
{\it Proof:} We write  $\, f=[1:f_1:\ldots :f_n]\,$ with holomorphic
$\,f_j: \bbbc \r \bbbc \se \{0\}$. Now we get with equations (\ref{1})
and (\ref{2}) for $j=1,\ldots ,n$:
$$
T(f_j,r) = T([1:f_j],r) + O(1) \leq T(f,r) + O(1),
$$
hence the $f_j$ are nonvanishing holomorphic functions of order
at most $\la$. This means that
$$
\lim\/{\rm sup}_{r \r \infty} \frac{T(f_j,r)}{r^{\la + \ep}} =0
$$
for any $\, \ep > 0$. From this equation our assertion follows with the
Weierstra\3 theorem as it is stated in \cite{HA}. \qed

The previous Lemma is helpful because we can use it to `calculate'
$T(f,r)$ by the Ahlfors-Lemma (cf.\ \cite{ST})
\begin{lem} \label{A}
Let $\, P_0,\ldots ,P_n\,$ be polynomials of degree at most $\,\la
\in \nz\,$.
Let $\, \al_j \in \bbbc$ be the coefficients of $\,x^{\lambda}\,$
in $P_j$ (possibly equal to zero). Let $\, L(\al_0,\ldots
,\al_n)\,$
be the length of the polygon defined by the convex hull of the
$\, \al_0,\ldots ,\al_n$. If
$$f=[e^{P_0}:\ldots :e^{P_n}]:\bbbc \r \pr_n$$
then
$$ \lim_{r \r \infty} \frac{T(f,r)}{r^{\la}} = \frac{L(\al_0,\ldots ,\al_n)}{2
\pi}$$
\end{lem}

We state the First and the Second Main Theorem of Value Distribution
Theory which relate the characteristic function and the counting
function (cf.\  \cite{SH}):

Let $\,f:\bbbc \r \pr_n\,$ be entire, and let $D$ be a divisor in $\pr_n$
of degree $d$, such that $\,f(\bbbc) \not\s \hbox{ {\rm support}}(D)$.
Then:

\medskip
{\bf First Main Theorem} $$ N_f(D,r) \leq d \cdot T(f,r) + O(1)$$

Another way of stating this theorem is the following:
The quantity
$$
\delta_f(D)=\liminf_{r\to \infty}\left(1- {N_f(D,r)\over d \cdot
T(f,r)}\right)
$$
is called {\it defect} of $D$ with respect to $f$. Then
$$
\delta_f(D)\geq 0.
$$

Assume now that $\, f(\bbbc)\,$ is not contained in any hyperplane in
$\pr_n$, and let $\, H_1,\ldots ,H_q\,$ be distinct hyperplanes in
general position.  Then

\medskip
{\bf Second Main Theorem}
$$
(q-n-1)T(f,r) \leq \sum_{j=0}^q N_f(H_j,r) + S(r)
$$
where $\: S(r) \leq O(\log (rT(f,r)))\,$ for all $\,r \geq r_0\,$ except
for a set of finite Lebesque measure. If $f$ is of finite order, then
$\, S(r) \leq O(\log r)\,$ for all $\,r \geq r_0$.

We examine how the characteristic function behaves under morphisms of the
projective space:
\begin{lem} \label{m}
Let
$$ R=[R_0:\ldots :R_N]: \pr_n \r \pr_N$$
be a morphism with components
of degree $p$, and let $\, f:\bbbc \r \pr_n\,$ be entire. Then
$$ T(R \circ f,r) = p \cdot T(f,r) + O(1)$$
\end{lem}
{\it Proof:} Define
$$\mu([z_0:\ldots :z_n]) = \frac{|R_0|^2+\ldots + |R_N|^2}{(|z_0|^2+\ldots
+|z_n|^2)^p}
$$
Since $R$ is a morphism the $\,R_j, j=0,\ldots ,N\,$ have no common zeroes,
hence there exist constants $\,A,B >0\,$ with
$\:0 < A \leq \mu \leq B\:$ on $\pr_n$. From that and equation (\ref{1})
we get:
$$T(R \circ f,r) - p \cdot T(f,r) = \frac{1}{4 \pi} \int_0^{2 \pi} (\log ||R
\circ
f||^2(re^{i\v}) - p \cdot \log ||f||^2(re^{i\v}))d \v +O(1) $$ $$=
\frac{1}{4 \pi} \int_0^{2 \pi} \log  (\mu \circ f)(re^{i\v}) d\v +O(1)$$
In the last term the integral is bounded by $\, \frac{1}{2}\log A\,$ and
$\,\frac{1}{2} \log B\,$
independently of $r$. \qed

\section{Some consequences of Brody's techniques}
In this section we list briefly some consequences of Brody's techniques
for later application.  The first is a corollary of a well known theorem
of M. Green.  It shows how to use entire curves $\,f:\bbbc \r \pr_2\,$ of
finite order to prove hyperbolicity of quasiprojective varieties.  The
second follows from of a theorem of M.Zaidenberg.

\medskip
a) The main theorem of \cite{GRE2} implies:
\begin{cor} \label{c}
Let $D$ be a union of curves $\, D_1,\ldots ,D_m\,$ in $\pr_2$ such that
for all $\, i=1,\ldots ,m\,$ the number of intersection points
of $\,D_i\,$ with $\:\bigcup_{j=1,\ldots ,m;j \neq i} D_j\:$ is at
least three.
Then $\,\pr_2 \se D\,$ is complete hyperbolic and hyperbolically
embedded, if there does not exist a non-constant entire curve
$\, f:\bbbc \r \pr_2\,$ of order at most two which misses $D$.
\end{cor}

b) The following proposition shows that the property of a union of curves
having
hyperbolic complement is essentially a (classically) open condition.
\begin{prop} \label{z}

Let $\,H_1,\ldots ,H_m\,$ be hypersurfaces in $\, \pr_2 \times (\D_t)^n\,$
for some $\,t>0$, $n \in \nz$. Let $\: \pi : \pr_2 \times (\D_t)^n
\r (\D_t)^n\:$ be the projection.  Assume that

1) for all $\, z \in (\D_t)^n\,$ and all $i=1,\ldots ,m$ the fibers
$\: \pi^{-1}(z) \cap H_i\:$ are curves in $\pr_2$

2) for all $i=1,\ldots ,m$ the number of intersection points of $\:
\pi^{-1}(0) \cap H_i\:$ and \\ $\:  \bigcup_{j=1,\ldots ,m;j \neq
i}(\pi^{-1}(0) \cap H_j))\:$ is at least three.

3) $\:  \pr_2 \se \bigcup_{j=1,\ldots ,m} (\pi^{-1}(0) \cap H_j)\:$ is
hyperbolically embedded in $\pr_2$.

Then $\:  \pr_2 \se \bigcup_{j=1,\ldots ,m}( \pi^{-1}(z) \cap H_j)\:$ is
complete hyperbolic and
hyperbolically embedded for all $\, z \in (\D_s)^n\,$ for some $\,s \leq
t \,$.
\end{prop}
{\it Proof:} In the terminology of \cite{ZA}, the $\: \pi^{-1}(0)
\cap H_i\:$ form an absorbing $H$-stratification (cf.\ \cite{ZA}, p. 354 f.),
for which we can apply Theorem~2.1 of \cite{ZA}.  Complete hyperbolicity
follows from \cite{LA2}, p.36. \qed

\section{Nonexistence of algebraic entire curves in generic
complements}
In this section we prove that the complement of 3 generic quadrics,
or of any 4 generic curves other than 4 lines, does not contain
non-constant entire curves contained in an algebraic curve.
 Because of Corollary \ref{c} this can
be regarded as a statement of `algebraic' hyperbolicity.

Let us first make precise what we mean by generic.  The space of curves
$ \Ga_i $ of degree $ d_i $ in $\pr_2$, which we define as the
projectivized space of homogeneous polynomials of degree $d_i$, is a
projective space of dimension $\:  n_i = \frac{1}{2}(d_i+2)(d_i+1) -1 $.
Hence $ {\cal C}(d_1,\ldots ,d_k) = \prod_{i=1}^k \pr_{n_i} $ is
projective algebraic.  In order to simplify notations we denote this
space by $S$ in all what follows, and its elements by $ s \in S$, and by
$ \Ga_i(s) $ the curve given by the i-th component of $ s \in S$.

\begin{prop} \label{a}
Let $ S={\cal C}(2,2,2) $ or $ S= {\cal C}(d_1,\ldots ,d_k) $ with $ k
\geq 4 $ and $ d= \sum_{i=1}^k d_i \geq 5 $. Then there exists a proper
algebraic variety $ V \s S $ st. for $ s \in S \se V $ the following
holds:\\ For any irreducible plane algebraic curve $A\subset \pr_2$ the
punctured Riemann surface $A\setminus \bigcup_{i=1}^k \Ga_i(s)$ is
hyperbolic, in particular any holomorphic map $f:\bbbc \to
\pr_2\setminus\bigcup_{i=1}^k \Ga_i(s)$ with $f(\bbbc)\subset A$ (which
may also be reducible) is constant.
\end{prop}
{\it Proof:} In order to define $  V \s S $ we list 5 conditions:

(1) All $ \Ga_i (s) $ are smooth (and of multiplicity one).

(2) The $ \Ga_i (s),\:  i=1,\ldots ,k $ intersect transversally,
in particular no  3 of these intersect in one point.

(3) In the case of $ {\cal C}(2,2,2) $:  For any common tangent
line of two of the quadrics $\Ga_j(s)$ which is tangential to these in
points $P$ and
$Q$ resp. the third quadric does not intersects the tangent in both points $P$
and
$Q$.

(4) In the case of $ {\cal C}(d_1,d_2,1,1),\:d_1,d_2 \geq 2 $:  There
does not exist a common tangent $L$ to $ \Ga_1(s) $ and $ \Ga_2(s) $
such that $L\cap\Ga_1(s)=\{P\}$ and $L\cap\Ga_2(s)=\{Q\}$ such that
the lines $ \Ga_3(s) $ and $\Ga_4(s) $ contain $P$ and $Q$ resp..

(5) In the case of $ {\cal C}(d_1,1,1,1) $:  There does not exist a
tangent line $L$ at $ \Ga_1 (s) $ with $L\cap \Ga_1 (s)= \{P\}$
such that $P$ is contained in one of the lines  $ \Ga_i (s),\:i=2,3,4 $
and $L$ contains the intersection points of the other two lines.

Define $ V \s S $ to be the set of those points $s \in S$ such that the
$\Ga_i(s)$ violates one of the above conditions.  This set is clearly
algebraic and not dense in $S$.

For intersections of at least five curves (2) implies that any
irreducible algebraic curve $A$ intersects $\bigcup_{i=1}^k \Ga_i(s)$ in at
least three
different points, which proves the claim.

Assume that there exists an irreducible algebraic curve $A \s \pr_2 $ and $s\in
S$
such that $A\setminus \bigcup_{i=1}^k \Ga_i(s)$ is not hyperbolic.
By condition (2) we know that $  A \cap \bigcup_{i=1}^k \Ga_i(s) $ consists of
at least
2 points $P$ and $Q$. Moreover, $A$ cannot have a singularity at $P$
or $Q$ with different tangents, because $A$ had to be reducible in such
a point, and $A\setminus \bigcup_{i=1}^k \Ga_i(s)$ could be identified with an
irreducible curve with
at least three punctures. (This follows from blowing up such a point or
considering the normalization).

So $ A \cap \bigcup_{i=1}^k \Ga_i(s) $ consists of exactly 2 points $P$ and $Q$
with simple tangents.  We denote the multiplicities of $A$ in $P$ and $Q$ by $
m_P $ and
$ m_Q $. Let $d_0= \deg(A)$.  Then the inequality (cf.\  \cite{FU}, p.117)
$$
m_P(m_P-1)+m_Q(m_Q-1) \leq (d_0-1)(d_0-2)
$$
implies
\begin{equation} \label{*}
m_P , m_Q < d_0 \hbox{ {\rm or} }d_0=m_P=m_Q=1.
\end{equation}
Let us now first treat the case $ k=4 $: Each $ \Ga_i(s) $ contains
exactly one of the points $P$ and $Q$. Let $\Ga_j(s)$ and $\Ga_k(s)$ resp.
intersect $A$ in $P$ and $Q$ resp. not tangential, i.e. with
tangents different from those of $A$ in these points. Let $d_j$ and $d_k$ be
the degrees of these
components. We compute intersection multiplicities according to \cite{FU},
p.75
$$ m_P = I(P, A \cap \Ga_j (s)) = d_j d_0\quad {\rm and} \quad m_Q=
I(Q, A \cap \Ga_k (s)) = d_k d_0.$$
Hence
$$ d_0=d_j=m_P=1\:\: {\rm and} \:\: d_0=d_k=m_Q=1.$$
In particular $A, \Ga_j$ and $\Ga_k$ are lines. These situations are
excluded by (4) and (5).

Now let us treat the case of 3 quadrics.  After a suitable
enumeration of its components we may assume that $P \in \Ga_1 (s)
\cap \Ga_2 (s) $ and $ Q \in \Ga_3 (s)$. If $Q \not\in \Ga_2 (s)
\cup \Ga_1 (s) $ we are done, since then we may assume that $A$ is
not tangential to $ \Ga_2 (s)$, and again
$$
m_P = I(P,A \cap \Ga_2 (s)) = 2d_0.
$$
contradicts equation (\ref{*}).
So we may assume that $\, Q \in \Ga_2(s) \cap \Ga_3(s)$. Now $A$ has to
be tangential to $\, \Ga_1(s)$ in $P$ and to $\, \Ga_3(s)$ in $Q$,
otherwise we again get $\, m_P=2d_0\,$ or $\,m_Q =2d_0\,$ what
contradicts equation (\ref{*}). But then $\, \Ga_2(s)$ is not tangential
to $A$ in $P$ or $Q$, so we have
$$m_P + m_Q = I(P,A \cap \Ga_2(s)) + I(Q, A \cap \Ga_2(s)) = 2d_0$$
Again by equation (\ref{*}) this is only possible if
$\, m_P=m_Q=d_0=1$, but then we are in a situation which we excluded
in condition (3), which is a contradiction. \qed

\section{Hyperbolicity of generic complements of at
least four curves}

In this section we prove a result in the direction towards a generalized
second main theorem. As a corollary we get a new proof of the fact that
for any generic collection of four hypersurfaces $\Gamma_j$,
$j=1,\ldots 4$ in $\pr_2$ any entire curve $f:\bbbc \to \pr_2 \setminus
\bigcup_{j=1}^4 \Gamma_j$ has to be algebraically degenerate. This fact,
combined with our result in the previous section implies the
hyperbolicity of the complement of such a configuration.

\begin{theo} \label{4c}
Let $  S= {\cal C} (d_1,\ldots ,d_k) $ with $  k \geq 4$, $d=
\sum_{i=1}^k d_i \geq 5$.  Then there exists an algebraic variety $  V
\s S $ such that for $ s \in S \se V $ the following holds:  Assume that
$ f:  \bbbc \r \pr_2 \se  \bigcup_{j=1}^3 \Ga_j (s)  $ is a
non-constant holomorphic map. Then $ \de_f (\Ga_l
(s))=0 $ for $l=4,\ldots ,k$.  In particular, $f$ cannot miss any
$\Ga_l(s)$, $l=4,\ldots ,k$.
\end{theo}
{\it Proof:} Let $ V \s S $ be defined like in the proof of Proposition
\ref{a}, i.e.  $ s \in S \se V $, iff the conditions (1) to (5) given
there are satisfied.   Let $ \Ga_i(s) = \{P_i(s) = 0 \} $ for
$i=1,\ldots ,k$. For suitable powers $a_j$
we have because of condition (2) a morphism
\begin{equation} \label{a1}
\Phi : \pr_2 \r \pr_2; [z_0:z_1:z_2] \r
[P_1^{a_1}(s):P_2^{a_2}(s):P_3^{a_3}(s)]
\end{equation}
Since there exists no non-constant morphism on projective spaces, whose
image is of lower dimension, for all $ s \in S \se V $ the image $ \Phi
(\pr_2) $ is not contained in an algebraic curve.  From now on, we keep
some $s \in S\se V $ fixed and drop the parameter $s$ for the rest of
the proof. Furthermore let $ \Phi (\Ga_4) =\{ Q=0 \} $, where

$$Q(w_0,w_1,w_2) = \sum_{i_0 +i_1+ i_2 =e} a_{i_0i_1i_2}
w_0^{i_0}w_1^{i_1}w_2^{i_2},$$
so $\deg Q=e$.  Finally, let
$\: \Phi^{-1}(\Phi (\Ga_4)) = \Ga_4 \cdot R \:$ be the decomposition
of the inverse image curve of the curve $ \Phi(\Ga_4) $ in $\Ga_4$
and the other components (which possibly may contain $\Ga_4$ as well).
Now the proof consists of 3 steps:\\

a) We have $ a_{e00} \not= 0,\: a_{0e0} \not= 0,\:a_{00e} \not= 0 $, i.e.
the polynomial $Q$ contains the $e$-th powers of the coordinates:\\
We prove that indirectly, so without loss of generality we may assume
that $ a_{e00}=0$.  Then we have $ Q([1:0:0])=0$, i.e.  $ [1:0:0] \in
\Phi (\Ga_4)$.  So there exists a point $ z \in \Ga_4 $ with $
P_1(z) \not= 0,\:P_2(z)=0,\:P_3(z)=0 $. But that means that the 3 curves
$\Ga_2$, $\Ga_3$ and $\Ga_4$ have a common point which contradicts our
condition (2).\\

b) We show by using the Second Main Theorem that $ \de_{\Phi \circ f}(
\Phi (\Ga_4))=0$:\\
Let $ J=\{ (\underline{i} = (i_0,i_1,i_2) :  a_{i_0i_1i_2} \not= 0 \}$
and $ \kappa :  J \r \{0,1,\ldots ,p \} $ be an enu\-me\-ra\-tion of
$J$.  Let $ Q_j = w_0^{i_0}w_1^{i_1}w_2^{i_2} $ if $ \kappa
((i_0,i_1,i_2)) =j$.  Then by part a) the map
$$
\Psi : \pr_2 \r \pr_p; [w_0:w_1:w_2] \r [Q_0:\ldots :Q_p]
$$
is a morphism with components of degree $ e=\deg(Q)$.  The $p+2$ lines
$ L_i= \{ \xi_i=0
\},\:i=0,\ldots ,p $ and $ L= \{\sum_{\underline{i} \in J} a_{\underline{i}}
\xi^{\kappa (\underline{i})} =0 \} $ are in general position.
Furthermore the map $\:  \Psi \circ \Phi \circ f :  \bbbc \r \pr_p $ is
linearly non degenerate:  By Proposition \ref{a}, $ f(\bbbc) $ is not
contained in an algebraic curve, so especially not in an algebraic curve
of the form $ \sum_{\underline{i} \in J} b_{\underline{i}}(P_1^{a_1})^{i_0}
(P_2^{a_2})^{i_1}(P_3^{a_3})^{i_2} $, resulting from such a line in $ \pr_p$,
unless the
latter is identically zero.  But this is impossible, since the map
$\Phi$ is surjective.  So we have by the Second Main Theorem:
$$T(\Psi \circ \Phi \circ f,r) \leq N_{\Psi \circ \Phi \circ f}(L,r) +
\sum_{i=0}^p N_{\Psi \circ \Phi \circ f}(L_i,r) + S(r) $$
and by the First Main Theorem
$$N_{\Psi \circ \Phi \circ f}(L,r) \leq T(\Psi \circ \Phi \circ f,r) +
O(1)$$
Observe that all $N_{\Psi \circ \Phi \circ f}(L_i,r)$ vanish.
Together with Lemma \ref{m} this yields, since $\deg Q=e$
\begin{equation} \label{a2}
\de_{\Phi \circ f} (\Phi (\Ga_4)) = \liminf_{r \r \infty} (1- \frac{N_{\Phi
\circ f}(\Phi(\Ga_4),r)}{\deg(Q)T(\Phi \circ f,r)}) =
\liminf_{r \r \infty} (1- \frac{N_{\Psi \circ \phi \circ f} (L,r)}{
T(\Psi \circ \Phi \circ f,r)}) =0
\end{equation}

c) We finally show that $ \de_f(\Ga_4)=0$:\\ By equation (\ref{a2}) and
Lemma \ref{m} we have:  $$1 = \limsup_{r \r \infty} \frac{N_{\Phi \circ
f}( \Phi (\Ga_4),r)}{ \deg(Q)T(\Phi \circ f,r)} = \limsup_{r \r \infty}
\frac{N_f(\Phi^{-1} \Phi (\Ga_4),r)}{\deg(Q \circ \Phi)  T(f,r)}$$ $$=
\limsup_{r \r \infty} \frac{N_f(\Ga_4,r) + N_f(R,r)}{\deg(Q \circ \Phi)
T(f,r)} = \limsup_{r \r \infty} \frac{N_f(\Ga_4,r) +
N_f(R,r)}{(\deg(\Ga_4) + \deg(R)) T(f,r)}$$ or short:
\begin{equation} \label{a3}
\limsup_{r \r \infty} \frac{N_f(\Ga_4,r)}{T(f,r)} + \limsup_{r \r \infty}
\frac{N_f(R,r)}{T(f,r)} = \deg(\Ga_4) + \deg(R)
\end{equation}
By the First Main Theorem we have:
$$ \limsup_{r \r \infty} \frac{N_f(R,r)}{\deg(R)T(f,r)} \leq 1 ,\:\:
\limsup_{r \r \infty} \frac{N_f (\Ga_4,r)}{\deg(\Ga_4)T(f,r)} \leq 1$$
and hence with equation (\ref{a3}):
$$\limsup_{r \r \infty} \frac{N_f(\Ga_4,r)}{T(f,r)} = \deg (\Ga_4)
,\:{\rm i.e.}\: \de_f(\Ga_4)=0 .$$
\qed

\section{Line systems through intersection points of three quadrics}

In this section, we study certain configurations of 18 lines associated to
three smooth quadrics.  These lines are needed in order to apply Value
Distribution Theory to prove our main theorem in the next section.

Let $V' \s S={\cal C}(2,2,2)$ be the algebraic variety defined by
the conditions (1), (2), and (3) given in the Proof of Proposition
\ref{a}, namely $s\in S \se V'$, iff

\begin{description}
\item[(1)] All $\Ga_i(s)$ are smooth quadrics.

\item[(2)] The $\Ga_i(s),\:i=1,2,3$ intersect transversally (in
particular not all 3 intersect in one point)

\item[(3)] For any common tangent line of two of the quadrics $\Ga_j(s)$
which is tangential to these in points $P$ and $Q$ resp. the third quadric
does not intersects the tangent in $P$ and $Q$.
\end{description}

In order to prove our main theorem we will need one further condition of
`genericity' related to those 18 lines already mentioned above.  For
this condition it is quite not so obvious any more that it yields a
quasiprojective set.  We shall give an argument for this in Proposition
\ref{ls}.

Let us first state the extra condition:  Because of (2) any two of the
three quadrics $\Ga_1$, $\Ga_2$, $\Ga_3$ intersect in 4 distinct points
$A_1, A_2, A_3, A_4$ which give rise to six lines \begin{equation}
\label{!} \o{A_1A_2},\,\o{A_3A_4} \hbox{ \rm and
}\o{A_1A_3},\,\o{A_2A_4} \hbox{ \rm and } \o{ A_1A_4},\,\o{A_2A_3}.
\end{equation} So all three pairs of quadrics give rise to three sets
$L_{12}(s)$, $L_{13}(s)$ and $L_{23}(s)$ of six lines each, i.e. a
collection $L(s)$ of 18 lines. We will show in the proof of Proposition
\ref{ls} that as a consequence of (1) and (2) they are pairwise distinct.

Now our condition (4) reads:

\begin{description}
\item[(4)]
The 18 lines  $L(s)$ intersect as follows: At
any point of $\Ga_i(s)\cap \Ga_j(s)$, $i \not= j$
there intersect exactly 3 of the 18 lines, and in every other point of
$\pr_2$ there intersect at most 2 of the 18 lines.
\end{description}

Now we have:
\begin{prop} \label{ls}
Define $V \s S$ to be the set of all $s \in S$ such that one of the
conditions (1) to (4) is not satisfied.  Then $V \s S$ is a
proper algebraic subset.
\end{prop}
{\it Proof:}\/
In order to prove the Proposition we use an argument which involves an
elementary case of a Chow scheme.

We denote by $\pr_2^\vee$ the space of all lines in $\pr_2$.  Look at
the following rational map
$$
\psi: (\pr_2)^4 \to (\pr_2^\vee)^6
$$
$$
(A_1,A_2,A_3,A_4) \mapsto (A_j \wedge A_k)_{j<k}
$$
where the wedge product of two points is considered as an element of the
dual projective space. This map descends to a rational map of symmetric
spaces:
$$
\Psi: S^4(\pr_2) \to S^6(\pr_2^\vee).
$$
Over the complement of a proper algebraic subset it assigns to a set of
four distinct points the configuration of six lines through these
points.

Now we assign to any $s \in S\se V'$ the tripel of sets
$(\Ga_1(s)\cap\Ga_2(s),\Ga_1(s)\cap\Ga_3(s),\Ga_2(s)\cap\Ga_3(s))$,
which amounts to a morphism
$$
\rho:S\se V' \to (S^4(\pr_2))^3.
$$
Observe that $\Xi:=(\Psi)^3\circ \rho : S\se V' \to
(S^6(\pr_2^\vee))^3$ is a morphism. Now we can rephrase condition (4):

Let $U\simeq \bbbc^3$ and $W\simeq (U)^6$.  Then we consider $W^3 = \{
(a_{jk})|a_{jk} \in U ; j=1,\ldots,6; k=1,2,3 \}$ and look at the linear
subspace $B\s W^3$ defined by the condition that at least three
components $a_{j_1 k_1}$, $a_{j_2 k_2}$ and $a_{j_3 k_3}$ are {\it
linearly dependent} where {\it not all $k_j$ are the same}. (We needn't
care about the system of the
six lines given by the four intersection points of
two fixed quadrics, since they automatically have the desired intersection
properties, because no three of the four intersection points of the
two quadrics can be collinear.) Obviously $B$
descends to an algebraic set $\tilde B \in (S^6(\pr_2^\vee))^3$.  Now
(4) means for $s\in S\se V'$ that $\Xi(s) \not \in \tilde B$.  The
construction immediately implies that $\,V \se V' \s S \se V'\,$ is
algebraic, and since $\, V' \s S$ is algebraic, we have that $\, V=
\overline{V \se V'}
\cup V'\,$ is algebraic in $S$, where the closure here means the Zariski
closure.

We have to show that $V \neq S$.  The existence of an $s \in S\setminus
V$ is proved by a deformation argument:  We start with any $s \in S \se
V'$ (then the $\Ga_i(s)$ are smooth and we have 12 different
intersection points of two of the 3 quadrics each).  It is easy to see
that then we really have 18 different lines, otherwise 4 of the 12
intersection points of the 3 quadrics had to be contained in a line
(because no three quadrics pass through a line).  It follows from the
construction, that this line would intersect one of the quadrics
in 4 points, which is impossible.

Let $k=k(s)$ be the largest number of lines among the 18 lines
(determined by the parameter $s$) which run
through some point.  Let $\nu_k=\nu_k(s)$ be the number of points in
$\pr_2$ which are contained in $k(s)$ of the lines.
We will proceed now as follows:  We observe that $k$
lines running through a point is a closed condition with respect to the
classical topology of $S$.  That means that in a neighborhood $U$ of a
point $s_0 \in S$ we have $k(s)\leq k(s_0)$, and at least
$\nu_k(s)\leq\nu_k(s_0)$, if $k(s)= k(s_0)$.  We will show that
for some $s\in U$ actually $k(s)<k(s_0)$ or at least
$\nu_k(s)<\nu_k(s_0)$, if $k(s)= k(s_0)$,  as long as $k(s_0)>3$
or $k(s_0)=3$ but $\nu_k(s_0)> 12 $.

Iterating this procedure we are done if we can show:  Consider the 18
lines in $L(s_0)$.  If $k \geq 4$ take any of these intersection points
where $k$ lines intersect (call it $T$), if $k=3$ take such an
intersection point $T$ which is not intersection point of two of the
quadrics.  Then we can find $s \in S$ arbitrarily near to $s_0$ st. over
$s$ the point $T$ `breaks up' into intersection points of strictly less
then $k(s_0)$ lines. But then $k(s)<k(s_0)$, or at least $k(s)=k(s_0)$
and $\nu_k(s)<\nu_k(s_0)$.

Let us now prove that:  Take 3 of the lines running through $T$ over
$s_0$ and denote them by $L_1,L_2,L_3$.  Each of them is defined by
construction by two of the intersection points of two of the 3 quadrics.
Let $L_1$ be defined by such points $T_1,T_2$, let $L_2$ be defined by
$T_3,T_4$ and let $L_3$ be defined by $T_5,T_6$.  We may assume that no
3 of the 6 points $T_1,\ldots ,T_6$ are equal to $T$ (this could only
occur if $T$ is an intersection point of 2 of our 3 quadrics, but then
$k \geq 4$ and we just have taken the 3 lines defined by $T$ and one
other intersection point each, so we can choose a different line).  So
without loss of generality we may assume that $T_1 \not= T \not= T_2$,
and we have the following 3 possibilities for
$L_1,L_2,L_3,T,T_1,\ldots ,T_6$:

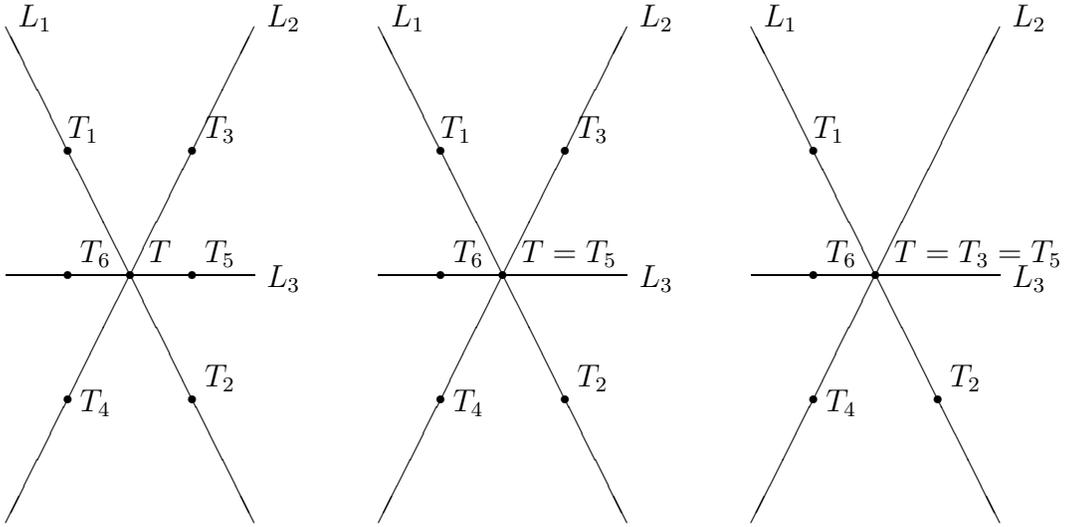
\moveleft1cm
\hbox{

    \beginpicture
     \setcoordinatesystem units <0.2em,0.2em>
     \unitlength0.2em

     \setplotarea x from -30 to 30, y from -50 to 50

     \put {\line(1,-2){20}} [Bl] at 0 0
     \put {\line(-1,-2){20}} [Bl] at 0 0
     \put {\line(1,0){20}} [Bl] at 0 0

     \put {\line(1,2){20}} [Bl] at 0 0
     \put {\line(-1,2){20}} [Bl] at 0 0
     \put {\line(-1,0){20}} [Bl] at 0 0

     \put {\circle*{1.1}} [Bl] at 0 0       
     \put {\circle*{1.1}} [Bl] at -10 20    
     \put {\circle*{1.1}} [Bl] at 10  -20   
     \put {\circle*{1.1}} [Bl] at 10 20     
     \put {\circle*{1.1}} [Bl] at -10 -20   
     \put {\circle*{1.1}} [Bl] at 10 0      
     \put {\circle*{1.1}} [Bl] at -10 0     

     \put {$L_1$} [Bl] at -18 40
     \put {$L_2$} [Bl] at  22 40
     \put {$L_3$} [Bl] at  22 -2

     \put{$T_1$} [Bl] at -10 22
     \put{$T_2$} [Bl] at  12 -18
     \put{$T_3$} [Bl] at  12 22
     \put{$T_4$} [Bl] at  -8 -22
     \put{$T_5$} [Bl] at  12  2
     \put{$T_6$} [Bl] at  -8  2
     \put{$T$}   [Bl] at   3 2

    \endpicture

    \beginpicture
     \setcoordinatesystem units <0.2em,0.2em>
     \unitlength0.2em

     \setplotarea x from -30 to 30, y from -50 to 50

     \put {\line(1,-2){20}} [Bl] at 0 0
     \put {\line(-1,-2){20}} [Bl] at 0 0
     \put {\line(1,0){20}} [Bl] at 0 0

     \put {\line(1,2){20}} [Bl] at 0 0
     \put {\line(-1,2){20}} [Bl] at 0 0
     \put {\line(-1,0){20}} [Bl] at 0 0

     \put {\circle*{1.1}} [Bl] at 0 0       
     \put {\circle*{1.1}} [Bl] at -10 20    
     \put {\circle*{1.1}} [Bl] at 10  -20   
     \put {\circle*{1.1}} [Bl] at 10 20     
     \put {\circle*{1.1}} [Bl] at -10 -20   
     \put {\circle*{1.1}} [Bl] at -10 0     

     \put {$L_1$} [Bl] at -18 40
     \put {$L_2$} [Bl] at  22 40
     \put {$L_3$} [Bl] at  22 -2

     \put{$T_1$} [Bl] at -10 22
     \put{$T_2$} [Bl] at  12 -18
     \put{$T_3$} [Bl] at  12 22
     \put{$T_4$} [Bl] at  -8 -22
     \put{$T_6$} [Bl] at  -8  2
     \put{$T=T_5$}   [Bl] at   3 2

    \endpicture

    \beginpicture
     \setcoordinatesystem units <0.2em,0.2em>
     \unitlength0.2em

     \setplotarea x from -30 to 30, y from -50 to 50

     \put {\line(1,-2){20}} [Bl] at 0 0
     \put {\line(-1,-2){20}} [Bl] at 0 0
     \put {\line(1,0){20}} [Bl] at 0 0

     \put {\line(1,2){20}} [Bl] at 0 0
     \put {\line(-1,2){20}} [Bl] at 0 0
     \put {\line(-1,0){20}} [Bl] at 0 0

     \put {\circle*{1.1}} [Bl] at 0 0       
     \put {\circle*{1.1}} [Bl] at -10 20    
     \put {\circle*{1.1}} [Bl] at 10  -20   
     \put {\circle*{1.1}} [Bl] at -10 -20   
     \put {\circle*{1.1}} [Bl] at -10 0     

     \put {$L_1$} [Bl] at -18 40
     \put {$L_2$} [Bl] at  22 40
     \put {$L_3$} [Bl] at  22 -2

     \put{$T_1$} [Bl] at -10 22
     \put{$T_2$} [Bl] at  12 -18
     \put{$T_4$} [Bl] at  -8 -22
     \put{$T_6$} [Bl] at  -8  2
     \put{$T=T_3=T_5$}   [Bl] at   3 2

    \endpicture
     }

\noindent
The point $T_1$ lies on 2 of the $\Ga_i$, assume on $\Ga_1 \cap \Ga_2$.
Then at most 3 of the 4 or 5 different points in $\{T_2,\ldots,T_6\}$
can also be in $\Ga_1 \cap \Ga_2$.  So there exists one of them, call it
$T_0$, which does lie on $\Ga_1 \cap \Ga_3$ or $\Ga_2 \cap \Ga_3$,
assume on $ \Ga_2 \cap \Ga_3$.  So at most 4 of the points $T_2,\ldots
,T_6$ are contained in $\Ga_1$.  So we can `move' $\Ga_1$ while keeping
these 4 points fixed and keeping $\Ga_2$ and $\Ga_3$ fixed.  But that
means that there is a non-constant variation of $s$ where we keep all of
the points $T_2,\ldots,T_6$ fixed. Hence the lines
$L_2$ and $L_3$ and their intersection point $T$ are kept fixed.  We
claim that for some small such variation the line $L_1$ does
not pass any longer through $T$.  If it would, it had to be fixed, since
$T_2$ is kept fixed.  By definition we have $T_1 \in L_1 \cap \Ga_1 \cap
\Ga_2$ and $ L_1 \cap \Ga_2$ is a discrete set.  Hence $T_1$ remains
fixed.  But that would mean that any quadric $\Ga_1$ through the at most
4 of the fixed points
$T_2,\ldots,T_6$ contained in $\Ga_1$ must contain a
fifth fixed point $T_1$.  This is certainly a contradiction, since the space
of plane quadrics is of dimension five.  \qed

\pagebreak
{}From any of the configurations of 18 lines in Proposition \ref{ls}
we can pick 12 in general position:

\begin{cor} \label{12l} There exists an algebraic variety $V
\s S$ st. for all $ s \in S \se V$ we have subsets of 12 of the 18 lines
of Proposition \ref{ls} which are in general position.
\end{cor}
{\it Proof:} For each pair $\Ga_i(s),\Ga_j(s),\:  i \not=j$ we have
constructed 3 pairs of lines (defined by equation (\ref{!})).  Choose,
for fixed $s \in S$, for each pair of quadrics two of these pairs of
lines.  \qed

At last we prove the simple fact that the pairs of lines as
defined in equation (\ref{!}) are contained in the linear system
spanned by the two quadrics.

\begin{prop} \label{sy}
Let $\Ga_1,\Ga_2$ be two smooth quadrics intersecting in 4 different
points $ A_1,A_2,A_3,A_4$, and let the lines $L_1$ resp.  $L_2$ be given
by $A_1,A_2$ and $A_3,A_4$ resp.  Then $L_1L_2$ is a degenerate
quadric contained in the linear system spanned by $\Ga_1$ and $\Ga_2$,
i.e.
$\:  L_1L_2 = a \Ga_1 + b \Ga_2$.  \end{prop}
{\it Proof:} Look at the set ${\cal L}$ of all quadrics (possibly
singular) which run through the 4 points $A_1,A_2,A_3,A_4$.  Then ${\cal
L}$ is a one dimensional linear system containing $L_1L_2$.  Since it is
one dimensional, it is spanned by any 2 of its elements, e.g.\ by $\Ga_1$
and $\Ga_2$. \qed

\section{Hyperbolicity of generic complements of three quadrics}
We will prove: \begin{theo} \label{mt} Let $V \s S$ be the variety
defined in Proposition \ref{ls}. Let $s \in S \se V$. Then the
quasiprojective variety $\: \pr_2 \se \bigcup_{i=1}^3 \Ga_i(s)\:$ is
complete hyperbolic and hyperbolically embedded.
\end{theo}

\begin{rem}
{\it The variety $S\setminus V$ is certainly not contained in an open
subset of the space of all divisors of degree $6$ whose complement in
$\pr_2$ is hyperbolic (cf.\ also \cite{ZA}):  Take any quadratic polynomials
$P_1,P_2, P_3$ corresponding to some $s\in S\setminus V$.
Then with respect to suitable coordinates we have $P_1=z_0^2-z_1z_2$.  Set
$Q=(z_1^6 + P_1\cdot F)$, where $F$ is an arbitrary polynomial of degree
$4$, $P=P_1\cdot P_2 \cdot P_3$, and $R_t= P + t\cdot Q$, $t\in \bbbc$.
Then the zero set of $R_0$ is just $\bigcup_{i=1}^3 \Ga_i(s)$.  However, for
$t\not=0$
the intersection of $V(R_t)$ with the rational curve $V(P_1)$ consists
only of the point $[0:0:1]\in \pr_2$.}
\end{rem} \qed

{\it Proof of the Theorem:} By Corollary \ref{c} it is sufficient to
show that there doesn't exist a non-constant entire curve $\:f:  \bbbc \r
\pr_2 \se \bigcup_{i=1}^3 \Ga_i(s)$ of order at most 2.

Assume there exists such a non-constant entire curve $f$.  From
Proposition \ref{a} we know that $f$ is not algebraically degenerate.

For simplicity of notation we drop the $s$ in the rest of the proof.
Furthermore we enumerate the 12 lines which we constructed in
Proposition \ref{ls} and Corollary \ref{12l} as follows:

$L_1L_2$ and $L_3L_4$ are in the linear system of $\Ga_1$ and $\Ga_2$

$L_5L_6$ and $L_7L_8$ are in the linear system of $\Ga_1$ and $\Ga_3$

$L_9L_{10}$ and $L_{11}L_{12}$ are in the linear system of $\Ga_2$ and
$\Ga_3$.

Let $\Ga_i =  \{P_i = 0\}$ with a homogeneous polynomial $P_i$
of degree 2.

The map $\: \Phi = [P_1:P_2:P_3]: \pr_2 \r \pr_2$ is a morphism
(because $\Ga_1 \cap \Ga_2 \cap \Ga_3 = \emptyset$). Furthermore the map
$\: \Phi \circ f : \bbbc \r \pr_2$ again is an entire curve and the
map $\Phi \circ f$ is again of finite order at most 2, because
by Lemma \ref{m} we have
\begin{equation} \label{*1}
T( \Phi \circ f,r) = 2 \cdot  T(f,r) +O(1)
\end{equation}
Since $f$ misses the divisor $\Ga_1\Ga_2\Ga_3$ the map $\Phi \circ
f$ misses the divisors $\{ z_i = 0 \}, \:i=1,2,3$ and hence
by Lemma \ref{e} we can write
\begin{equation} \label{*2}
\Phi \circ f = [g_0:g_1:g_2]
\end{equation}
with
$$
g_i = e^{\al_i \xi^2 + \be_i \xi + \ga_i};\: \al_i,\:\be_i,\:\ga_i
\in \bbbc
$$
where $ \:g_i = (P_i \circ f) \cdot  h$; $h: \bbbc \r \bbbc^*$
are entire functions.

We may assume that not all three $\al_j$ are equal:  Assume $\al_1 =
\al_2 = \al_3$, then we can divide out the function $e^{\al_1 \xi^2}$
and then compose the resulting functions with $\xi\mapsto \xi^2$, i.e.
we may consider the function $ \Phi \circ f (\xi^2)$.  This map is again
of order at most 2 and we have $\:  g_i = e^{\be_i \xi^2 + \ga_i}$.  If
now $\be_1 = \be_2 = \be_3$, the map $ \Phi \circ f$ would be constant,
which is impossible, since $f$ is algebraically non degenerate.
So we exclude the case $\al_1=\al_2 = \al_3$   without loss of
generality.

The Ahlfors Lemma \ref{A} allows the  computation of some limits of
characteristic functions: For $1 \leq i < j \leq 3$:
\begin{equation} \label{*3}
\lim_{r \r \infty} \frac{T([P_i \circ f: P_j \circ f],r)}{r^2}
= \lim_{r \r \infty} \frac{T([g_i:g_j],r)}{r^2}
= \frac{2|\al_i - \al_j|}{2 \pi}
\end{equation}
and
\begin{equation} \label{*4}
\lim_{r \r \infty} \frac{T([P_1 \circ f: P_2 \circ f: P_3 \circ f],r)}{
r^2} = \lim_{r \r \infty} \frac{T([g_1:g_2:g_3])}{r^2}
\end{equation}
$$
= \frac{|\al_1 - \al_2| + |\al_1 - \al_3| + |\al_2 - \al_3|}{2 \pi}
$$
hold.  Now we want to relate the counting functions of the 12 lines to
the characteristic functions used in equations (\ref{*3}) and
(\ref{*4}):  We know that $L_1L_2$ is in the linear system of $\Ga_1$
and $\Ga_2$, i.e.  $\:  L_1L_2 = a \Ga_1 + b \Ga_2\:$ with $a,b \not= 0$
since $\Ga_1$ and $\Ga_2$ are smooth quadrics. We consider the map
$$[P_1 \circ f: P_2 \circ f] : \bbbc \r \pr_1.$$
Its image is not contained in a hyperplane in $\pr_1$, i.e.
a point, since $f$ is algebraically non degenerate. Furthermore the 3
divisors $$\: [z_0 = 0], [z_1 = 0], [az_0 + bz_1 = 0]$$ are in
general position in $\pr_1$,  i.e. distinct. The  Second Main Theorem
yields
$$T([P_1 \circ f: P_2 \circ f],r) \leq
N_{[P_1 \circ f: P_2 \circ f]} ([z_0 = 0],r)$$ $$ +
N_{[P_1 \circ f: P_2 \circ f]} ([z_1 = 0],r)
+ N_{[P_1 \circ f: P_2\circ f]} (az_0 + bz_1 =0],r) + O(\log r) =$$
$$
N_f ([P_1 = 0],r) + N_f([P_2=0],r) + N_f([aP_1 + bP_2 =0],r) + O(\log r)=$$
$$0+0+N_f([L_1L_2=0],r) +O(\log r) = N_f([L_1=0],r) + N_f([L_2=0],r) +
O(\log r)$$
where $N_f([P_i=0],r)=0$ because $f$ misses $ \Ga_i = [P_i=0]$,
and where we identify the line $L_i$ with its defining equation, so that
$[L_i = 0]$ makes sense. On the other hand we have by the First Main
Theorem
$$N_{[P_1 \circ f: P_2 \circ f]}([az_0 + bz_1 =0],r) \leq
T([P_1 \circ f: P_2 \circ f],r) + O(1)$$
and hence
$$T([P_1 \circ f: P_2 \circ f],r) = N_f([L_1=0],r) + N_f([L_2=0],r) + O(\log r)
$$ The corresponding equations hold for all other lines as well, i.e. we
have:
$$T([P_1 \circ f:P_2 \circ f],r) = N_f([L_1=0],r) + N_f([L_2=0],r)+ O(\log r)$$
                              $$   = N_f([L_3=0],r) + N_f([L_4=0],r)+ O(\log r)
$$
\begin{equation} \label{*5}
  T([P_1 \circ f:P_3 \circ f],r) = N_f([L_5=0],r) + N_f([L_6=0],r)+ O(\log r)
\end{equation}
              $$                   = N_f([L_7=0],r) + N_f([L_8=0],r)+ O(\log
r)$$
$$T([P_2 \circ f: P_3 \circ f],r)= N_f([L_9=0],r) + N_f([L_{10}=0],r) +O(\log
r)$$
                         $$      =N_f([L_{11}=0],r) + N_f([L_{12}=0],r) +O(\log
r).
$$
Since $ f: \bbbc \r \pr_2$ is not linearly degenerate and the 12 lines
$L_1,\ldots ,L_{12}$ are in general position, we can again apply the
Second Main Theorem and get
\begin{equation} \label{*6}
9 \cdot T(f,r) \leq \sum_{i=1}^{12} N_f([L_i=0],r) + O(\log r).
\end{equation}
The equations (\ref{*1}), (\ref{*5}) and (\ref{*6}) imply
$$ \frac{9}{2} \cdot T(\Phi \circ f,r) = 9 \cdot T(f,r) +O(1) \leq
\sum_{i=1}^{12}
N_f([L_i=0],r) + O(\log r)$$
$$= 2 \cdot (T([P_1 \circ f:P_2 \circ f],r) + T([P_1 \circ f:P_3 \circ f],r)
+ T([P_2 \circ f: P_3 \circ f],r) + O(\log r)$$
Hence together we have
\begin{equation} \label{*7}
9 \cdot T(\Phi \circ f,r) \leq 4 \cdot ( \sum_{1 \leq i < j \leq 3} T([P_i
\circ f:
P_j \circ f],r)) + O(\log r).
\end{equation}
We now divide equation (\ref{*7}) by $r^2$ and take
$\lim_{r \r \infty}$. Using the equations (\ref{*3}) and
(\ref{*4}) we obtain:
$$ 9 \cdot \frac{|\al_1 - \al_2| + |\al_1 - \al_3| + |\al_2 - \al_3|}{2 \pi}
\leq 4 \cdot 2 \cdot \frac{|\al_1 - \al_2| + |\al_1 - \al_3|+|\al_2 - \al_3|}{2
\pi}.
$$
This can only hold if $\al_1 = \al_2 =\al_3 $,  which is a contradiction.
\qed

\def\sq{\hbox{\rlap{$\sqcap$}$\sqcup$}}

\def\qed{\ifmmode\sq\else{\unskip\nobreak\hfil
\penalty50\hskip1em\null\nobreak\hfil\sq
\parfillskip=0pt\finalhyphendemerits=0\endgraf}\fi}

\def\bbbr{{\rm I\!R}} 
\def\bbbn{{\rm I\!N}} 
\def\bbbm{{\rm I\!M}}
\def\bbbh{{\rm I\!H}}
\def\bbbk{{\rm I\!K}}
\def\bbbp{{\rm I\!P}}

\def\bbbc{{\mathchoice {\setbox0=\hbox{$\displaystyle\rm C$}\hbox{\hbox
to0pt{\kern0.4\wd0\vrule height0.9\ht0\hss}\box0}}
{\setbox0=\hbox{$\textstyle\rm C$}\hbox{\hbox
to0pt{\kern0.4\wd0\vrule height0.9\ht0\hss}\box0}}
{\setbox0=\hbox{$\scriptstyle\rm C$}\hbox{\hbox
to0pt{\kern0.4\wd0\vrule height0.9\ht0\hss}\box0}}
{\setbox0=\hbox{$\scriptscriptstyle\rm C$}\hbox{\hbox
to0pt{\kern0.4\wd0\vrule height0.9\ht0\hss}\box0}}}}

\def\bbbz{{\mathchoice {\hbox{$\sf\textstyle Z\kern-0.4em Z$}}
{\hbox{$\sf\textstyle Z\kern-0.4em Z$}}
{\hbox{$\sf\scriptstyle Z\kern-0.3em Z$}}
{\hbox{$\sf\scriptscriptstyle Z\kern-0.2em Z$}}}}

\def\bbbq{{\mathchoice {\setbox0=\hbox{$\displaystyle\rm Q$}\hbox{\raise
0.15\ht0\hbox to0pt{\kern0.4\wd0\vrule height0.8\ht0\hss}\box0}}
{\setbox0=\hbox{$\textstyle\rm Q$}\hbox{\raise
0.15\ht0\hbox to0pt{\kern0.4\wd0\vrule height0.8\ht0\hss}\box0}}
{\setbox0=\hbox{$\scriptstyle\rm Q$}\hbox{\raise
0.15\ht0\hbox to0pt{\kern0.4\wd0\vrule height0.7\ht0\hss}\box0}}
{\setbox0=\hbox{$\scriptscriptstyle\rm Q$}\hbox{\raise
.15\ht0\hbox to0pt{\kern0.4\wd0\vrule height0.7\ht0\hss}\box0}}}}

\def\blk{\hbox{ \vrule height 7pt width 4pt depth 0pt} }

\newtheorem{expl}[defi]{Example}




\section{Complements of two quadrics and a line}

In this section we need the following theorem of M.~Green \cite{GRE3}
(in degree $d=2$) which generalizes in a sense the classical Borel
lemma.

\begin{theo}\label{grebor}
a) Let $g_0$, $g_1$, $g_2$ be entire holomorphic functions of finite
order, $g_1$ and $g_2$ both nowhere vanishing. Assume that
$$ g_0^2+ g_1^2+g_2^2 = 1. \eqno{(A)}$$
Then the set
$$\{1,g_0,g_1,g_2\}$$
of holomorphic functions has to be linearly dependent.\\
b) Let $g_0$ and $g_1$ be entire holomorphic functions of finite order,
$g_1$ nowhere vanishing. Assume that
$$ g_0^2+ g_1^2 = 1. $$
Then $g_0$ and $g_1$ must be constant.
\end{theo}

We consider the complement of three quadrics.  We allow one of these to
be also a double line.  (The case, where one of the three quadrics
degenerates to two distinct lines, i.e. two quadrics and two lines,
has already been treated above).

Since in this section we work also with double lines we will distinguish
between $\Ga$ and $P$, where $\Ga=V(P)$ (for simplicity reasons we
didn't always do this in the previous sections).


Before we state the main result of this section, we observe that
also some singular configurations of two quadrics in the projective
plane can be treated by means of the generalized Borel lemma.


\begin{prop}\label{propQ}
Let $\Gamma_j=\{Q_j=0\}\subset \bbbp_2$, $j=1,2$ be two smooth distinct
quadrics, whose intersection consists of exactly one point.  Then any
holomorphic map $f:\bbbc \to \bbbp_2\setminus (\Gamma_1\cup \Gamma_2)$
of finite order has values in a quadric (which may degenerate to a
double line) from the linear system spanned by $Q_1, Q_2$.
\end{prop}
{\it Proof.}\/ Let the common tangent line to $\Gamma_1$ and $\Gamma_2$
through the intersection point
be defined by the linear equation $L=0$. One verifies immediately that
\begin{equation}
L^2=aQ_1+bQ_2,
\end{equation}
 $a,b \not = 0$. Let $q_j$ be entire non-vanishing
holomorphic functions, $j=1,2$ such that $q_j^2 = Q_j \circ f$, and
$q_0=L\circ f$. Then Theorem~8.1 b) implies that $Q_1 \circ f = c \cdot
Q_2 \circ f$.
\qed

Another case is the following.

\begin{prop}
Let $\Gamma_j=\{Q_j=0\}\subset \bbbp_2$, $j=1,2$ be two smooth distinct
quadrics, which intersect exactly at two points tangentially.  Then any
holomorphic map $f:\bbbc \to \bbbp_2\setminus (\Gamma_1\cup \Gamma_2)$
of finite order has values in a quadric contained in the linear system
spanned by $Q_1, Q_2$.
\end{prop}
{\it Proof.}\/ The linear system spanned by $Q_1$ and $Q_2$
contains $L^2$,
where $L$ is the line through the two points of intersection. Using
this the statement follows as above. \qed


\begin{theo}\label{borhyp}
Let $0\not=Q_j \in \bbbc[z_0,z_1,z_2]$, $j=1,2,3$ be quadratic
polynomials, where either all $Q_j$ are irreducible or all but one
which may be a square of a linear function.
Let $\Gamma_j\subset \bbbp_2$ be the zero-sets. Assume

\begin{description}
\item[(1)] no more than two of these intersect at one point,

\item[(2)] no tangent to a smooth quadric $\Gamma_j$ at a point of
intersection with some other $\Gamma_k$
contains a further intersection point of the curves $\Gamma_l$,

\item[(3)] there exists a linear combination of the $Q_j$ which is a
square:
\begin{equation}\label{Sum}
\sum_{j=1}^3 a_j Q_j = P^2,\quad   P\in \bbbc[z_0,z_1,z_2],
\end{equation}
where at least two coefficients $a_j$ are different from zero.
\end{description}

Then any holomorphic map
$$
f:\bbbc \to \bbbp_2\setminus \bigcup_{j=1}^3 \Gamma_j
$$
has values in a quadric (which may be degenerate to a double line).
\end{theo}

We call a holomorphic map $\bbbc \to \bbbp_2$ {\it linearly or
quadratically degenerate}\/, if its values are contained in a line or a
(possibly degenerate) quadric resp..

\begin{cor}\label{corhyp}
\phantom{abc}
Let $\Gamma_j=V(Q_j)\subset \bbbp_2$, $j=2,3$ be smooth quadrics and
$\Gamma_1=L_1=V(Q_1)\subset \bbbp_2$ a line, where $Q_1$ is the square
of a linear polynomial, and let the assumptions of \ref{borhyp} be satisfied.

\begin{description}
\item[(1)]
The quasiprojective variety $\bbbp_2\setminus
\bigcup_{j=1}^3 \Gamma_j$ is Brody-hyperbolic,
unless there exists a smooth quadric or a line $\Gamma$  such that
after choosing the notation accordingly ($p$, $q$ distinct points):

\begin{description}
\item[(a)]
$\Gamma \cap \Gamma_2=\{p,q\}$, $\Gamma \cap\Gamma _3= \{p\}$,  $\Gamma
\cap L_1=\{q\}$

\item[(b)]
$\Gamma \cap \Gamma_2=\{p\}$, $\Gamma \cap\Gamma _3= \{p\}$,  $\Gamma
\cap L_1=\{q\}$

\item[(c)]
$\Gamma \cap \Gamma_2=\{p\}$, $\Gamma \cap\Gamma _3= \{q\}$,  $\Gamma
\cap L_1=\{p\}$

\item[(d)]
$\Gamma \cap \Gamma_2=\{p\}$, $\Gamma \cap\Gamma _3= \{q\}$,  $\Gamma
\cap L_1=\{p,q\}$
\end{description}

\item[(2)]
The quasiprojective variety $\bbbp_2\setminus \bigcup_{j=1}^3 \Gamma_j$
is complete hyperbolic and hyperbolically embedded, unless

\begin{description}
\item[(e)]
at least two of the $\Gamma_j$ are tangent to each other at some point,

\item[(f)]
there exists a smooth quadric, which has only one point of intersection
with each of $\Gamma_2$ and $\Gamma_3$ with both of these points contained in
$\Gamma_1$,

\item[(g)]
There exists a  tangent  to one of $\Gamma_2$ and $\Gamma_3$ at a point of
intersection with $\Gamma_1$ which is tangent to the other smooth quadric.
\end{description}

\end{description}

\end{cor}  \qed


We introduce the following polynomials which will take care of a
necessary elimination process in the proof of \ref{borhyp}.

\begin{defi}\label{defR}
Let the homogeneous polynomial $R_j(y_0,\ldots,y_j)\in
\bbbc[y_0,\ldots,y_j]$ of degree $2^{j-1}$ be defined by the equation
$$R_j(x_0^2,\ldots,x_j^2)=
\prod_{(\epsilon_1,\ldots,\epsilon_j)\in \{1,-1\}^j}(x_0+ \epsilon_1 x_1
+ \ldots + \epsilon_j x_j).$$
\end{defi}

For later applications we need some properties of the $R_j$:

\begin{lem}\label{lemR}
a) $\; R_2 (x,y,z)= x^2 +y^2 +z^2 -2xy -2xz -2yz$\\
b) Let $a,b,c \in \bbbc$. Then
$$
S(x,y,z):=R_3(ax+by+cz,x,y,z)
$$
has the following properties:

\begin{description}
\item[1)]
The coefficient of $x^4$ equals $(a-1)^4$.
\item[2)]
The coefficient of $x^2y^2$ equals
$2(3a^2(b-1)^2-2a(b-1)(3b+1)+3b^2+2b+3)$. In particular,
if the coefficient of $y^4$ vanishes, the coefficient of
$x^2y^2$ equals $16$.
\item[3)]
Assume that all coefficients of forth powers in $S$ vanish. Then
$$
S(x,y,z)= 16(x^2y^2+x^2z^2+y^2z^2)-32(x^2yz+xy^2z+xyz^2).
$$
\end{description}
\end{lem}
We omit the computational proof. \qed
{\it Proof of Theorem \ref{borhyp}.}
Since the map $f$ has no values in the given quadrics $\Gamma_j$, there
exist entire holomorphic functions $q_j$, $j=1,2,3$ such that $q_j^2=
Q_j\circ f$.  If we put then $q_0=P\circ f$ and $g_j=q_j/q_3$ for
$j=0,1,2 $. We apply the generalized Borel lemma (Theorem 8.1): If one
of the $a_j$ vanishes,  we get immediately quadratic degeneracy from
part b) of this theorem.

So from now on we assume that all $a_j$ are non-zero.  Thereom~8.1~a)
implies that the set of functions $\{1,g_0,g_1,g_2\}$ is linearly
dependent, i.e.  $\{q_0,\ldots,q_3\}$ has this property.

Let
\begin{equation}\label{Sum1}
\sum_{j=0}^3 \alpha_j q_j =0, \quad {\hbox{not all }} \alpha_j=0,
\end{equation}
and let $R=R_3$ be
the polynomial of \ref{defR}.  It has been chosen in a way such that
$R(\alpha_0^2q_0^2,\ldots,\alpha_3^2q_3^2)=0$. The assumption
(\ref{Sum})  means
that $q_0^2=a_1q_1^2+a_2q_2^2+a_3q_3^2$. Now the curve defined by the
equation
\begin{equation}\label{RGl}
\tilde R(z_0,z_1,z_2)=
R(\alpha_0^2(a_1Q_1+a_2Q_2+a_3Q_3),\alpha_1^2Q_1,\alpha_2^2Q_2,
\alpha_3^2Q_3)=0
\end{equation}
contains the image of $f$ and is of degree at most eight. We have
to show that $\tilde R$ is not identically zero. Otherwise, since
$(Q_1,Q_2,Q_3)$ defines a morphism, i.e. an epimorphism
$Q:\bbbp_2\to \bbbp_2$,
the polynomial $R(\alpha_0^2(a_1 y_1^2 + a_2 y_2^2 + a_3 y_3^2), \alpha_1^2
y_1^2,
\alpha_2^2 y_2^2,\alpha_3^2 y_3^2) \in \bbbc[y_1,y_2,y_3]$ would be the
zero polynomial. The definition of $R$ would imply that
$\alpha_0^2(a_1 y_1^2 + a_2 y_2^2 + a_3 y_3^2) =
(\sum_1^3 \delta_j \alpha_j y_j)^2$ for certain $\delta_j=\pm1$.
Thus at least two of $\alpha_1, \alpha_2, \alpha_3$ must vanish.
However, by assumption, the $a_j$ are different from
zero. From this fact it follows immediately that all $\alpha_j=0$,
which is a contradiction. We have shown that $f(\bbbc)$ is contained
in an algebraic curve of degree at most eight which is defined by a
polynomial of degree four in $Q_1,Q_2,Q_3$. \qed

Before we proceed with the proof of Theorem~\ref{borhyp}, we give an
application
of the classical Borel Lemma.  Let $\bbbp_2=\{z_0,z_1,z_2\}$, and
$H_j=\{z_j=0\}$ be the coordinate hyperplanes.

\begin{rem}\label{Rem1}
Let $f:\bbbc \to \bbbp_2\setminus (H_0\cup H_1\cup H_2)$ be a
holomorphic map. Assume that $f$ is algebraically degenerate, i.e. its
values are contained in an algebraic curve $C$. Then $f(\bbbc)\subset C'$,
where $C'$ is the zero-set of a polynomial of the form $z_0^k-\beta
z_1^lz_2^m$, $\beta\not =0$, $k,l+m \leq deg(C)$
 (after a suitable reordering of indices).
\end{rem}
{\it Proof.}\/ Let $C$ be the zero-set of some
polynomial $P(z_0,z_1,z_2)$.  Denote by $f_j$ the components of $f$.
The classical lemma of Borel, applied to the monomials in the expansion
of $P$,
implies that there exist at least two such monomials, which are
proportional after
composing with $f$. Thus $f_0^r f_1^s f_2^t = \beta f_0^u f_1^v f_2^w$
for some $\beta\neq 0$.
\qed

(The statement of the Lemma has an obvious generalization to $\bbbp_n$.)

An immediate consequence of Theorem \ref{borhyp}, as far as we have proved is
yet,
and of Remark \ref{Rem1}  is, since all $Q_j\circ f$
have no zeroes,

\begin{lem}\label{lem1}
Given the assumptions of {\rm\ref{borhyp}}, the image $f(\bbbc)$ is contained
in
a curve of the form
$$Q_u^k-\alpha Q_v^l Q_w^m = 0, \quad \alpha\not = 0,$$
where $\{u,v,w\} =\{1,2,3\}$, and $k, l+m \leq 4$.
\end{lem} \qed

We note the following fact:

\begin{lem}\label{lem2}
Let $f:\bbbc \to \bbbp_2$ be as above.
\begin{description}
\item[(1)]
Let $f(\bbbc)$ be contained in the zero-set
\begin{equation}\label{equiv}
V(Q_1^{k_1}Q_2^{k_2}Q_3^{k_3}- \alpha Q_1^{l_1}Q_2^{l_2}Q_3^{l_3}),
\end{equation}
with $\sum k_j = \sum l_j = 4$,
and $k_j=l_j$ for at least one $j$,

or
\item[(2)]
let $f(\bbbc)$ be contained in both zero-sets
\begin{equation}\label{Gl1}
V(Q_1^{k_1}Q_2^{k_2}Q_3^{k_3}- \alpha Q_1^{l_1}Q_2^{l_2}Q_3^{l_3})
\end{equation}
and
\begin{equation}\label{Gl2}
V(Q_1^{m_1}Q_2^{m_2}Q_3^{m_3}- \beta Q_1^{n_1}Q_2^{n_2}Q_3^{n_3})
\end{equation}
with $\sum k_\nu =\sum l_\nu=\sum m_\nu=\sum n_\nu=4$ such that the
vector
$$
(k_1-l_1,k_2-l_2,k_3-l_3)
$$
is not a rational multiple of the vector
$$
(m_1-n_1,m_2-n_2,m_3-n_3).
$$
\end{description}
Then $f(\bbbc)$ is contained in a quadric curve, which is a member of
the linear system generated by two of the quadrics $Q_j$.
\end{lem}

We call the monomials
$Q_1^{k_1}Q_2^{k_2}Q_3^{k_3}$ and $Q_1^{l_1}Q_2^{l_2}Q_3^{l_3}$
satisfying (\ref{equiv}) {\it equivalent with respect to $f$}.
{\it Proof.}\/ We can eliminate one of the $Q_j$, say $Q_1$, and obtain
that $f(\bbbc)$ is contained in $V(Q_2^r - \gamma Q_3^r)$ for some
integer $r$, and $\gamma \in \bbbc$, because the $Q_j\circ f$
have no zeroes.  Taking roots we find that $f(\bbbc)\subset V(\sigma Q_2
- \tau Q_3)$, $\sigma,\tau\in \bbbc$ not both equal to zero.  \qed

We return to the proof of Theorem~\ref{borhyp}, and consider under which
conditions the above lemma can be applied. Let the situation of
Theorem~\ref{borhyp} be given. Denote by
$$
T(Q_1,Q_2,Q_3)=
R(\alpha_0^2(a_1Q_1+a_2Q_2+a_3Q_3),\alpha_1^2Q_1,\alpha_2^2Q_2,
\alpha_3^2Q_3)
$$
the polynomial of (\ref{RGl}).

We know that $T$ is not the zero-polynomial.  We already reduced
the proof of Theorem~\ref{borhyp} to the case, where all $a_j$ in (\ref{Sum})
are
different from zero.

{\it First part:}\/ Let all $\alpha_j \neq 0$. Thus we can (after
normalizing these constants to $1$) apply Lemma~\ref{lem2}.

{\it First case:}\/ We claim that the conditions of Lemma \ref{lem2} (2)
are satisfied, if at least two of the coefficients of $Q_j^4$ in $T$,
say those of $Q_1^4$ and $Q_2^4$, are different from zero.

If $Q_1^4$ and $Q_2^4$ are equivalent, $Q_1\circ f$ is a constant
multiple of $Q_2\circ f$, and $f$ is quadratically degenerate.

Otherwise there exist exponents  $(r_1,r_2,r_3)$, $(s_1,s_2,s_3)$
of $Q_j$ that match $(4,0,0)$ and $(0,4,0)$ resp. in the sense
of (\ref{Gl1}) and (\ref{Gl2}) resp. because of the classical Borel
Lemma.
Assume that the assumptions of Lemma~\ref{lem2}~(2) are not fulfilled,
so there exists a rational number $c$ such that
$$
(4,0,0)-(r_1,r_2,r_3)= c((0,4,0)-(s_1,s_2,s_3)).
$$
Since $(r_1,r_2,r_3)\not=(4,0,0)$ we have $r_1<4$, thus $s_1>0$ and
$c<0$.  Now $0 \geq -r_3=c(-s_3)\geq 0$ implies $r_3=s_3=0$,
so we can apply Lemma \ref{lem2} (1).

{\it Second case:}\/ The next case to consider is, where exactly one
forth power occurs, say $Q_1^4$.  According to Lemma~\ref{lemR}~2), the
coefficient of $Q_1^2Q_2^2$ in $T$ must be different from zero.

Assume first that $Q_1^4$ is equivalent
to $Q_1^2Q_2^2$ with respect to $f$. Then
Lemma~\ref{lem2}~(1) is applicable. If these monomials are not
equivalent, we have some (non trivial) relations
$Q_1^4\sim Q_1^{r_1}Q_2^{r_2}Q_3^{r_3}$ and $Q_1^2Q_2^2\sim
Q_1^{s_1}Q_2^{s_2}Q_3^{s_3}$. If the assumptions of Lemma~\ref{lem2}~(2)
would not hold,  we had
$$
(4,0,0)-(r_1,r_2,r_3)= c((2,2,0)-(s_1,s_2,s_3))
$$
for some $0\neq c \in \bbbq$ and $0\leq r_j,s_j \leq 4$, $\sum r_j=\sum
s_j=4$.  For $r_3=0$ Lemma~\ref{lem2}~(1) could be applied.  Only
$0<r_3\leq 4$ is left; in particular $-r_3=c(-s_3)$ implies $c>0$, $s_3
> 0$.  Now $r_1\neq 4$.  Thus $4-r_1=c(2-s_1)$ gives $s_1=0$ or $s_1=1$.
Furthermore $-r_2= c(2-s_2)$ holds.  Again $r_2=0$ makes \ref{lem2}~(1)
applicable so that we are left with $s_2=3$ or $s_2=4$.  Thus
$(s_1,s_2,s_3)=(0,3,1)$.  Hence $f$ has values in the quartic curve
$Q_1^2-\gamma Q_2 Q_3=0$ for some $\gamma\in \bbbc $.

Let $C$ be the curve $V(Q_1^2- \gamma Q_2 Q_3)$. We note first that
$$
C\cap (\Gamma_1\cup\Gamma_2\cup\Gamma_3)=
\Gamma_1 \cap (\Gamma_2\cup\Gamma_3)
$$
The case, where two smooth quadrics $\Gamma_j$ intersect in exactly one
point, yields immediately quadratic degeneracy by
Proposition~\ref{propQ}, and we are done.  If one of the $\Gamma_j$ is a
line, it cannot be tangent to both of the further given smooth quadrics
--- this is also excluded by assumption (2).  Thus $C\cap
(\Gamma_1\cup\Gamma_2\cup\Gamma_3)$ consists of at least three points.
As $f(\bbbc)$ is contained in
$C\setminus(\Gamma_1\cup\Gamma_2\cup\Gamma_3)$, the curve $C$ cannot be
irreducible unless $f$ is constant.  We are left with the case where $C$
decomposes into a line $l$ and a cubic.  We have
$C\cap\Gamma_1=\Gamma_1\cap(\Gamma_2\cup\Gamma_3)=
C\cap(\Gamma_1\cap(\Gamma_2\cup\Gamma_3))$, which implies $l\cap
\Gamma_1=(l\cap \Gamma_2)\cup(l\cap \Gamma_3)$.
This equality means that $l\neq \Gamma_1$ and that $l \cap \Gamma_1$
consists of two distinct points $p'$ and $p''$, (since no more than two
of the $\Gamma_j$ pass through a point).  Let $l\cap\Gamma_2=\{p'\}$ and
$l\cap\Gamma_3=\{p''\}$.  This means that $l$ is at least tangent to one
of the smooth quadrics and passes through one further intersection point
of the $\Gamma_j$.  This was excluded by assumption (2).


{\it Third case:}\/ Assume finally that all coefficients of $Q_j^4$ in
$T$ vanish.  According to Lemma~\ref{lemR}~3) the non-zero monomials in
$T$ are
$Q_j^2Q_k^2$, $j\neq k$ and $Q_j^2Q_kQ_l$, where $(j,k,l)$ run through
all cyclic permutations of $(1,2,3)$.  We pick $Q_1^2Q_2^2$ and check to
which of the monomials it can be equivalent with respect to $f$.
Lemma~\ref{lem2}~(1) is directly applicable to all possible cases but
$Q_1^2Q_2^2 \sim Q_1Q_2Q_3^2$ which implies $Q_1Q_2 \sim Q_3^2$. This
case was already treated.

The claim is now shown under the assumption that all $\alpha_j$ are
different from zero.

If two or more of the $\alpha_j$ vanish, the claim is already clear from
(\ref{Sum1}): We then get the equation $\:\alpha_j q_j
= -\alpha_k q_k\,$, which, after squaring both sides,
yields us quadratic degeneracy immediately, or if $q_0$ is involved,
by using that at least two of the $a_i$ are not zero.
  The remaining case is, where exactly one $\alpha_j=0$.
Here we
use $R_2$ from \ref{defR} and arrive at a polynomial $U(y_1,y_2,y_3)$ of
degree two, such that $f(\bbbc)$ is contained in the zero-set of
$U(Q_1,Q_2,Q_3)$.  Again Borel's lemma is applied to its monomials.  A
non-empty subset of $\{Q_1^2,Q_2^2,Q_3^2,Q_1Q_2,Q_1Q_3,Q_2Q_3\}$ has to
be divided into sets of $f$-equivalent polynomials.  In the view of
\ref{lem2}~(1) the only
case to remain is $Q_j^2\sim Q_kQ_l$ where $(j,k,l)$ is a cyclic
permutation of $(1,2,3)$.  This case was treated above.
\qed

In the  sequel we treat the case of the complement of two plane quadrics
and a line and the case of three Fermat quadrics. We show that

\begin{theo}\label{thm221}
There exist
\begin{enumerate}
\item[(a)]
a quasiprojective set $V\subset {\cal C}(2,2,1)$ of codimension one and
\item[(b)]
an open, non-empty subset $U\subset {\cal C}(2,2,1)$ containing $V$
\end{enumerate}
such that for all $s\in U$ the space $\bbbp_2\setminus \Gamma(s)$
is complete hyperbolic and hyperbolically embedded.
\end{theo}
{\it Proof.} The set $V$ will be constructed in a such a way that
the configurations $\Gamma(s)$ for $s\in V$ satisfy the conditions of
Proposition~3.2 so that $(b)$ will follow from the first statement.

Let $\bbbp_2=\{[z_0:z_1:z_2]\}$ and

\begin{eqnarray}
l &=& c_0 z_0 + c_1 z_1 + c_2 z_2 \label{B0}\\
Q_0&=&l^2 \label{B1}\\
Q_j&=&\sum_{k=0}^2 a_{jk}z_k^2 + b_{j0} z_0z_1 + b_{j1} z_0z_2+ b_{j2}
z_1z_2\label{B2}
\end{eqnarray}
for $j=1,2$.

We include $Q_0$ in this notation and compute  $a_{0k}, b_{0k}$
in terms of $c_l$.
 We shall discuss, when (\ref{Sum}) holds for these.

Let $A=(a_{jk})$ and $B=(b_{jk})$.
Let $\hat{A}$ be the adjoint matrix of $A$, i.e. $\hat{A} \cdot A =
\hbox{det}(A) E$.
For $\kappa^2,\lambda^2,\mu^2 \in \bbbc$
we consider the following linear combination
\begin{equation}\label{B3}
(\kappa^2,\lambda^2,\mu^2)\cdot
\left( \hbox{det}(A)
\left(
\begin{array}{c}
z_0^2\\
z_1^2\\
z_2^2
\end{array}
\right) +
\hat{A} B
 \cdot
\left(
\begin{array}{c}
z_0z_1\\
z_0z_2\\
z_1z_2
\end{array}
\right)
\right) =
(\kappa^2,\lambda^2,\mu^2)\cdot
\hat{A}\cdot
\left(
\begin{array}{c}
Q_0\\
Q_1\\
Q_2
\end{array}
\right)
\end{equation}

Looking at the left hand side one verifies that this expression
is a square of a linear polynomial, if
and only if the following equation holds:

\begin{equation}\label{B4}
(\kappa^2,\lambda^2,\mu^2)\cdot \hat{A} B = 2\hbox{det}(A) (\kappa
\lambda,
\kappa
\mu,
\lambda\mu)
\end{equation}
We set $a=(a_{jk})_{j>0} \in \bbbc^6$, $b=(b_{jk})_{j>0} \in \bbbc^6$,
and $c=(c_l)\in \bbbc^3$. So $A$, $B$ and $\hat{A}$ are now given
in terms of $a,b,c$. We define $M\subset \pr_2 \times \bbbc^3
\times \bbbc^6 \times \bbbc^6$ to be the set of all points
$([\kappa:\lambda:\mu],c,a,b)$ for which (\ref{B4}) holds.

For all $m\in M$ the inequality ${\rm dim}_m M \geq 14$ holds, since
(\ref{B4}) consists of three equations in $\kappa, \lambda, \mu, A,\hat{A},B$
and hence in $\kappa, \lambda, \mu,  a, b, c$. Consider
the canonical projection ${\rm pr}: \bbbp_2 \times \bbbc^3 \times
\bbbc^6 \times \bbbc^6 \to \bbbc^3 \times \bbbc^6 \times \bbbc^6$. Let
$c_0=(1,0,0)$, $a_0=
\left(
\begin{array}{ccc}
0 & 1 & 0 \\
0 & 0 & 1 \\
\end{array}
\right)$,
and $b\in \bbbc^6$ arbitrary. Then we calculate that
$
(\bbbp_2\times \{(c_0,b,a_0)\})\cap M
$
is zero-dimensional.

In Example \ref{expl} we shall give an explicit example of  a point $m_0=
([\kappa_0:\lambda_0:\mu_0],c_0,b_0,a_0)$ which is contained in such a
zero dimensional set,  where
$\kappa_0,\lambda_0, \mu_0 \neq 0$. Denote by $M_0\subset M$ an irreducible
component of $M$ containing $m_0$.  Now ${\rm pr}(M_0)\subset \bbbc^3 \times
\bbbc^6 \times \bbbc^6$ is algebraic and at least of dimension 14, because
the fiber is zero dimensional. (One
can check easily that ${\rm pr}(M_0) \neq \bbbc^3 \times\bbbc^6
\times\bbbc^6$).

Let
$$
(\kappa^2,\lambda^2,\mu^2)\cdot \hat{A}=(\phi,\psi,\chi),
$$
and $N=V(\phi \cdot \psi \cdot \chi \cdot {\rm det}(A))\subset M$.
 Observe $M_0\setminus N
\neq \emptyset$, since $m_0 \not\in N$ (what can be checked easily).
 Let $V'\subset
\bbbc^3\times\bbbc^6\times\bbbc^6$ be the quasi projective hypersurface
$V'={\rm pr}(M_0)\setminus {\rm pr}(N) \subset {\rm pr}(M_0\setminus N)$,
which is not empty:  The fiber of ${\rm pr}|M_0$ at $m_0$ is of
dimension zero, hence ${\rm dim (pr} (N)) \leq {\rm dim} N < {\rm dim} M_0
= {\rm dim (pr}(M_0))$.

By means of the assignment $\bbbc^3\times\bbbc^6\times\bbbc^6 \ni
(c,b,a)\mapsto (l,Q_1,Q_2) \in \bbbc^3\times\bbbc^6\times\bbbc^6 $
we associate to any point of ${\rm pr}(M_0)$ a triple consisting of one linear
and two quadratical polynomials. Now  ${\rm pr}(M_0)$ as well as ${\rm pr}(N)$
are invariant under the canonical action of $(\bbbc^*)^3$, given by
multiplication of $l$, $Q_1$, $Q_2$ by elements of $\bbbc^*$.
This follows
from the original definition of $M$ and $N$ (the existence of a linear
combination of the $Q_0, Q_1, Q_2$ to a square and the number of
coefficients which are zero is independent of the $\bbbc^*$ action on
$l$, $Q_1$, $Q_2$) and the fact that under this action
$(\bbbc^*)^3\times M_0$ has values in some irreducible component of
$M$, which has to be $M_0$.

Now ${\rm pr}(M) \setminus {\rm pr}(N)$
defines a quasi projective subvariety $V' \subset {\cal C}(1,2,2)$ of
codimension one. Our aim is to construct a quasi projective variety $V
\subset {\cal C}(1,2,2)$ of codimension one, which is contained in
$V'$ satisfying the further conditions of Corollary~\ref{corhyp}(2), and
hence proving the Theorem. We already chose $M$ and $N$ in a way that
$V'$ satisfies condition (3) of \ref{borhyp}. All of the configurations
which had to be excluded because of the further conditions in \ref{borhyp} and
\ref{corhyp} define a proper
algebraic subset $W\subset {\cal C}(1,2,2)$. All we need is to see that
$V:= V'\setminus W$ is not empty. But we have ${\rm pr}(m_0) \in
V' \setminus W$ for our point $m_0$ coming from the example below.

\begin{expl}\label{expl}
The following set of quadratic polynomials defines an element of $V$. In
particular the complement of its zero-sets in $\bbbp_2$ is complete
hyperbolic and hyperbolically embedded.
\begin{eqnarray}
Q_0&=& z_0^2\\
Q_1&=&z_1^2+ z_0z_1 + z_0z_2 + (1/25) z_1z_2 \\
Q_2&=&z_2^2 + 50 z_0z_1 - 10 z_0z_2 + 9z_1z_2
\end{eqnarray}
\end{expl}

One checks immediately that $225 Q_0 + 100 Q_1 + 4 Q_2$ is a square.
Set $\Gamma_j=V(Q_j)$. Furthermore:

1) No more than two $\Gamma_j$ intersect in one point.

2) None of the $\Gamma_j$ are tangent to any other $\Gamma_k$.

3) No tangent to one of $\Gamma_2$ and $\Gamma_3$ at a point of
intersection with any $\Gamma_j$ contains a further point of
intersection of the $\Gamma_j$.

4) No tangent to one of $\Gamma_2$ and $\Gamma_3$ at a point of
intersection with $\Gamma_1$ is tangent to the other smooth quadric.

5) There exists no smooth quadric $\Gamma$ with $\, \Gamma_2 \cap \Gamma
= \{p'\}$, $\Gamma_3 \cap \Gamma = \{p''\}$ and $\{p',p''\} \s \Gamma_1$.\\

How to check 1) to 4) is obvious. If $T'$ resp. $T''$ are the linear
polynomials which give the  tangents
at $\Gamma_2$ in $p'$ resp. at $\Gamma_3$ in $p''$ we have
$\:Q = a Q_1 + b(T')^2\,$, $\: Q = c Q_3 + d (T'')^2\,$, where
$\Gamma = V(Q)$, $\Gamma_i = V(Q_i)$.
Now solve for $a,b,c,d$, and show that only the trivial solution
exists.     \qed

For intersections of three smooth quadrics Theorem~\ref{borhyp} is not quite
superseeded by the more general statement of Theorem~7.1. as the
application to intersections to Fermat quadrics shows. We first note a
further corollary to Theorem~\ref{borhyp}.

\begin{cor}

Let $\Gamma_j=V(Q_j)\subset \bbbp_2$, $j=1,2,3$ be smooth quadrics,
and let the assumptions of \ref{borhyp} be satisfied.

\begin{description}
\item[(1)]
The quasiprojective variety $\bbbp_2\setminus
\bigcup_{j=1}^3 \Gamma_j$ is Brody-hyperbolic,
unless there exists a smooth quadric or a line $\Gamma$  such that
after choosing the notation accordingly ($p$, $q$ distinct points):

\begin{description}
\item[(a)]
$\Gamma \cap \Gamma_1=\{p,q\}$, $\Gamma \cap\Gamma _2= \{p\}$,  $\Gamma
\cap \Gamma_3=\{q\}$

\item[(b)]
$\Gamma \cap \Gamma_1=\{p\}$, $\Gamma \cap\Gamma _2= \{p\}$,  $\Gamma
\cap \Gamma_3=\{q\}$
\end{description}

\item[(2)] The above conditions (a) and (b) can be replaced by the
following (somewhat stronger) condition:
\begin{description}
\item[(c)]
all of the $\Gamma_j$ intersect transversally.
\end{description}
In this case $\bbbp_2\setminus \bigcup_{j=1}^3 \Gamma_j$ is complete
hyperbolic and hyperbolically embedded.

\end{description}

\end{cor} \qed

We apply the Corollary to the following

\begin{prop}
Let
$$
Q_j=a_j x^2 + b_j y^2 + c_j z^2 ; \quad j=1,2,3
$$
be linearly independent polynomials, whose zero-sets $\Gamma_j$ are
smooth. Assume

\begin{description}
\item[(1)] no more than two of the $\Gamma_j$ intersect at one point,

\item[(2)] no tangent to a quadric $\Gamma_j$ at a point of
intersection with some other $\Gamma_k$
contains a further intersection point of the curves $\Gamma_l$,

\item[(3)]
none of the $\Gamma_j$ are tangent to each other at any point.
\end{description}
Then $\bbbp_2\setminus \bigcup_{j=1}^3 \Gamma_j$ is complete hyperbolic
and hyperbolically embedded.
\end{prop}
\qed


\pagebreak

\noindent Gerd Dethloff\\
Mathematisches Institut der Universit"at G"ottingen\\
Bunsenstra\3e 3-5\\
3400 G"ottingen\\
Germany\\
\vspace{0.8cm}e-mail: DETHLOFF@CFGAUSS.UNI-MATH.GWDG.DE\\

\noindent Georg Schumacher\\
Ruhr-Universit"at Bochum, Fakult"at f"ur Mathematik\\
Universit"atsstra\3e 150\\
4630 Bochum 1\\
Germany\\
\vspace{0.8cm}e-mail:GEORG.SCHUMACHER@RUBA.RZ.RUHR-UNI-BOCHUM.DE\\

\noindent Pit-Mann Wong\\
Dept. of Mathematics, University of Notre Dame\\
Notre Dame, Indiana 46556\\
USA\\
e-mail:PMWONG@CARTAN.MATH.ND.EDU

\end{document}